\documentclass[12pt,letterpaper]{article}
\usepackage{amssymb,amsfonts,amsthm,amsmath}

\usepackage{epsf}
\usepackage{graphicx}
\usepackage{placeins}

\def\hang{\hangindent\parindent}
 \def\rf{\par\noindent\hang}
 
\begin{document}

\baselineskip=20pt

\begin{center} \large{{\bf UPPER BOUNDS ON THE MINIMUM COVERAGE PROBABILITY OF
CONFIDENCE INTERVALS IN REGRESSION AFTER MODEL SELECTION}}
\end{center}

\bigskip

\begin{center}
\large{
{\bf PAUL KABAILA$^{1*}$ {\normalsize \rm AND} KHAGESWOR GIRI$^1$}}
\end{center}

\begin{center}
{\large{\sl La Trobe University}}
\end{center}

\vspace{11cm}

\noindent * Author to whom correspondence should be addressed.

\noindent $^1$ Department of Mathematics and Statistics, La Trobe University, Victoria 3086, Australia.

\noindent e-mail: P.Kabaila@latrobe.edu.au

\noindent Facsimile: 3 9479 2466

\noindent Telephone: 3 9479 2594


\newpage
\begin{center}
\large{{\bf Summary}}
\end{center}

\noindent We consider a linear regression model, with the parameter of interest a specified 
linear combination of the regression parameter vector. We suppose that, as a first step,
a data-based model selection (e.g. by preliminary hypothesis tests or minimizing AIC)
is used to select a model. It is
common statistical practice to then construct a confidence interval for the parameter
of interest based on the assumption that the selected model had been given to us {\it a priori}.
This assumption is false and it can lead to a confidence interval with poor coverage properties.
We provide an easily-computed finite sample
upper bound (calculated by repeated numerical evaluation of a double integral) 
to the minimum coverage probability of this confidence interval.
This bound applies for model selection by any of the following methods: minimum AIC,
minimum BIC, maximum adjusted $R^2$, minimum Mallows' $C_P$ and $t$-tests.
The importance of this upper bound is that it delineates general categories of 
design matrices and model selection procedures for which this confidence interval
has poor coverage properties. This upper bound is shown to be a finite sample analogue
of an earlier large sample upper bound due to Kabaila and Leeb.


\bigskip

\noindent {\sl Key words:} Adjusted $R^2$-statistic; AIC; 
``Best subset'' regression; BIC;  Mallows' criterion; $t$-tests.



\newpage

\begin{center}
{\large {\bf 1. Introduction}}
\end{center}


 It is very common in applied statistics that the model initially proposed is relatively complicated.
 The standard statistical methodology  
 for simplifying a complicated model is to carry out a preliminary data-based
 model selection by, for example, using preliminary hypothesis tests or minimizing AIC. This 
 is usually followed by the inference of interest, using the same data, based on the assumption
 that the selected model had been given to us {\it a priori}. This assumption is false and it
 can lead to an inaccurate and misleading inference. 
 In one particular context, Breiman (1992) has called this ``a quiet scandal in the 
 statistical community''.
 Nonetheless, this type of inference is taught extensively in university 
 courses and is applied widely in practice.
 It is therefore important to ascertain the extent to which this type of inference is 
 inaccurate and misleading.

 Consider the important case that the inference of interest is either a confidence interval
 or a confidence region. A confidence interval (region)
 with nominal coverage $1-\alpha$ that is constructed after preliminary model selection,
 using the same data and based on the (false) assumption that the selected model 
 had been given to us {\it a priori}, will be called a
 `naive' $1-\alpha$ confidence interval (region). 
 The literature on the coverage properties of naive confidence 
 intervals and regions is relatively recent. 
  Regal \& Hook (1991) provide an example 
 of a log-linear model, parameters
 and model selection procedure for which the coverage probability of the naive 
 0.95 confidence interval is far below 0.95. Hurvich \& Tsai (1990) provide examples
 of a linear regression model, parameters and model selection procedures for which the 
 naive 0.9, 0.95 and 0.99 confidence regions for the regression parameter vector have coverages far below 
 0.9, 0.95 and 0.99 respectively. 
 These authors do not seek to provide a comprehensive analysis of the coverage probability functions
 of the confidence intervals or regions they consider.
 Arabatzis et al. (1989),  Chiou \& Han (1995a, b),
 Chiou (1997) and Han (1998) find the minimum coverage probabilities of naive 
 confidence intervals in the contexts of some simple models and simple model selection procedures. 
 The minimum coverage probability of the naive confidence interval can be calculated for simple
model selection procedures in linear regression involving only a single variable
(Kabaila (1998)). The kinds of model selection procedures used
in practice in linear regression are typically much more complicated. For the real-life
example considered by Kabaila (2005), there are 20 variables each of
which is to be either included or not, leading to
a choice from among $2^{20}$
different models. 
In more complicated situations such as these,  
Kabaila (2005), Kabaila \& Leeb (2006, Section 3) and Giri \& Kabaila (2007)
use Monte Carlo simulation methods to assess the minimum coverage probability of 
the naive confidence interval, in the context of linear regression models.
A model selection procedure is said to be `consistent' if, for any fixed model parameters and 
sample size $\rightarrow \infty$,
 the true order of the model is consistently estimated.
 Minimization of BIC is such a procedure.
Kabaila (1995) and Leeb \& P\"otscher (2005)
are concerned with dispelling the misconception 
that naive confidence intervals and regions, constructed after a consistent preliminary model selection, 
 will have good coverage properties provided that the sample size is sufficiently large.

 Whilst this literature provides examples 
 of the poor coverage performance
 of naive confidence intervals, it may still be asked whether these examples are merely oddities
 or whether they are indicative of a more widespread phenomenon. 
 The way to answer this
 question is by delineating general categories of models
 and model selection procedures 
 for which the naive confidence interval has poor coverage properties. 
The aim of the present paper is to make a contribution to such a delineation in the 
context of the complicated type of model selection procedures used in practice for the linear regression model
 \begin{equation*}
 Y=X\beta + \varepsilon
 \end{equation*}
where $Y$ is a random $n$-vector of responses, $X$ is a known
$n\times p$ matrix with  linearly independent columns, $\beta$ is
an unknown parameter $p$-vector and $\varepsilon \sim N(0, \sigma^2 I_n)$
where $\sigma^2$ is an unknown positive parameter. Suppose that
the quantity of interest is $\theta =a^T \beta$ where $a$ is a
known $p$-vector ($a\neq 0$). Our aim is to find a confidence
interval for $\theta$ with minimum coverage probability a pre-specified value
$1-\alpha$, based on an observation of $Y$.

We suppose that, as a first step, a data-based model selection
is used to select a model. 
Specifically, suppose that the model selection procedure is used to
either set $\beta_i$ equal to 0 or allow it to vary freely for
each $i=q +1, \dots, p \; (q \geq 1)$.
We consider a confidence interval for $\theta$ with nominal
coverage $1-\alpha$ constructed on the (false) assumption that the
selected model had been given to us {\it a priori}. This is the
naive $1-\alpha$ confidence
interval for $\theta$.
Let $\hat \Theta$, $\hat \beta_{q+1}, \ldots, \hat \beta_p$ denote the least squares estimators
of $\theta$, $\beta_{q+1}, \ldots, \beta_p$ respectively.
Let Corr$(\hat \Theta, \hat \beta_j)$ denote the correlation between $\hat \Theta$ and $\hat \beta_j$.
Assume, without loss of generality, that $|\text{Corr}(\hat \Theta, \hat \beta_j)|$
is maximized with respect to
$j \in q +1, \dots, p$ at $j=p$. 
We use $\rho$ to denote the important parameter $\text{Corr}(\hat \Theta, \hat \beta_p)$.

We call a model selection procedure `conservative' when it is not consistent but, for any fixed model parameters,
the probability of choosing only correct models converges to 1 as the sample size $\rightarrow \infty$.
Kabaila \& Leeb (2006) provide an easily-computed large sample
upper bound (calculated by repeated numerical evaluation of a single integral) 
to the minimum coverage probability of this confidence interval
for conservative model selection procedures. Minimization of AIC is such a procedure.
Consider the case that a conservative model selection procedure
is used. 
The large sample upper bound of Kabaila \& Leeb (2006) is a 
continuous decreasing function of $|\rho|$, which approaches 
0 as $|\rho|$ approaches 1 from below.
This result tells
us is that for large samples, the naive $1-\alpha$ confidence interval
has minimum coverage probability far below $1-\alpha$ when $|\rho|$
is close to 1.
The importance of this result is that it delineates general categories of 
design matrices $X$ and model selection procedures
for which the naive confidence interval has poor coverage properties in large samples.

In the present paper we provide an easily-computed 
finite sample analogue (calculated by repeated numerical evaluation of a double integral) 
of the large sample upper bound
of Kabaila \& Leeb (2006). This finite sample upper bound applies to a wide range 
of model selection procedures, and is not restricted to conservative ones.
For conservative model selection procedures the large sample upper bound
complements the finite sample bound nicely.
We suppose that the model selection is based on one of the following methods:
(a) minimum AIC, (b) minimum BIC, (c) maximum adjusted $R^2$-statistic, (d) minimum Mallows' $C_P$
and (e) for each $j \in \{q+1, \ldots, p \}$ a $t$-test of the null hypothesis $H_{0j}: \beta_j=0$ 
against the alternative hypothesis $H_{Aj}: \beta_j \ne 0$.
We provide a method for obtaining an upper bound on the
minimum coverage probability of the naive confidence interval as follows.

For convenience, we introduce the following terminology. If the model selection 
procedure is (hypothetically) used to either set $\beta_i$ equal to 0 or allow it to vary freely for
each $i \in L$, where $L$ is a proper subset of $\{q+1,\ldots, p\}$, 
then we say that ``the model selection procedure is applied only to 
$\beta_i \in L$''.  
The following result is proved in section 2.
For each given $\ell$ satisfying $q < \ell < p$, the minimum coverage probability of the naive $1-\alpha$ confidence interval is bounded above by the
coverage probability of the naive $1-\alpha$ confidence interval for given
$\frac{1}{\sigma}(\beta_{\ell+1},\ldots, \beta_p)$
and the model selection procedure applied only to $\beta_{\ell+1},\ldots, \beta_p$.
Therefore, the minimum coverage probability of the
naive confidence interval is bounded above by the coverage probability of the naive
$1-\alpha$ confidence interval for given $\frac{1}{\sigma} \beta_p$ and the model selection
applied only to $\beta_p$. 
In Section 3 we derive an easily-computed expression for this upper bound for given $\frac{1}{\sigma} \beta_p$.
This expression is easily minimized numerically with respect to $\frac{1}{\sigma} \beta_p$ to obtain
the value of the finite sample upper bound on the minimum coverage probability of the naive $1-\alpha$ confidence interval.

This upper bound is a continuous decreasing function of $|\rho|$.
Some illustrative numerical evaluations of this upper bound are presented in Section 4.
See, for example, Figure 1 which is a plot of this upper bound as a function
of $|\rho|$ for model selection by minimizing Mallows' $C_P$,
with $m=n-p=5, 20, 50, 1000$ and $\infty$ (i.e. the large sample upper bound of Kabaila
\& Leeb (2006)).
The new finite sample upper bound tells us that the naive $1-\alpha$ confidence interval
has minimum coverage probability far below $1-\alpha$ when $|\rho|$
is close to 1.
The importance of this result is that it delineates a general category of 
design matrices $X$ and model selection procedures
for which the naive confidence interval has poor coverage properties in finite samples.


\bigskip

\begin{center}
{\large {\bf 2. Two important preliminary results}}
\end{center}

Suppose that the model selection procedure is used to
either set $\beta_i$ equal to 0 or allow it to vary freely for
each $i=q +1, \dots, p \; (q \geq 1)$. 
Let $\tilde{\cal K}$ denote the family of all subsets of $\{ q+1, \dots, p \}$,
including the empty set $\emptyset$.
We use $\tilde{K}$ to denote the element of $\tilde{\cal K}$
chosen by the model selection procedure.
Let $\hat{\beta}$ denote the least-squares estimator of $\beta$.
Let $RSS$ denote the following residual sum of squares,
$$RSS=(Y-X \hat{\beta})^T (Y-X\hat{\beta}).$$
Let $K$ be a fixed subset of $\{q +1, \dots, p\}$ and suppose
that $\beta_i$ is set equal to zero for each $i \in K$ and is
freely-varying for each $i\notin K$. Let $\vert K\vert$ denote the
number of elements in $K$. Also let $H_K$ denote the $\vert K
\vert \times p$ matrix whose $i$th row consists of zeros except
for the $j$th element which is 1 where $j$ is the $i$th ordered
element of $K$. Thus $H_K \beta =0$. Let $\hat{\beta}_K$ denote
the least-squares estimator of $\beta$ subject to this restriction.
Also let $RSS_K$ denote the residual sum of squares
$$RSS_K =(Y-X \hat{\beta}_K)^T(Y-X\hat{\beta}_K),$$
and $S_K^2 = RSS_K/(n-p+|K|)$.
The standard $1-\alpha$ confidence interval for $\theta$, assuming
that $H_K\beta =0$, is
\begin{equation*}
I(K) = \big[ a^T \hat{\beta}_K - d_K, a^T\hat{\beta}_K  + d_K \big]
\end{equation*}
where $d_K = t(n-p+\vert K\vert ) S_K \sqrt{v(K)}$, 
$t(m)$ is defined by $P \big (-t(m) \leq T\leq t(m) \big )=1-\alpha$ for
$T \sim {\rm t}_m$ and $v(K)$ is defined to be (variance of
$a^T\hat{\beta}_K)/\sigma^2$.

We consider the following 4 methods of model selection.

\medskip

\noindent \underbar{Method 1} (minimizing an AIC-like criterion)

\noindent 
$\tilde K$ minimizes
\begin{equation*}
AIC(K) =n \ln (RSS_K)+2 (p-\vert K\vert )f(n)
\end{equation*}
with respect to $K \in \tilde{\cal K}$. Here, $f(n)$ is 1 for AIC and $\frac{1}{2}\ln(n)$ for
BIC.

\medskip

\noindent \underbar{Method 2} (minimizing Mallows' $C_P$)

\noindent 
$\tilde K$ minimizes
\begin{equation*}
C_K=\frac{RSS_K}{RSS/(n-p)} -n+2(p-\vert K\vert )
\end{equation*}
with respect to $K \in \tilde{\cal K}$.

\medskip

\noindent \underbar{Method 3} (maximizing adjusted $R^2$)

\noindent 
$\tilde K$ minimizes
\begin{equation*}
B_K=\frac{RSS_K}{n-p+|K|} 
\end{equation*}
with respect to $K \in \tilde{\cal K}$. 

\medskip

\noindent \underbar{Method 4} ($t$-tests)

\noindent 
$\tilde K$ consists of the set of $j \in \{q+1, \ldots, p \}$ for which a $t$-test
of the null hypothesis $H_{0j}: \beta_j=0$ 
against the alternative hypothesis $H_{Aj}: \beta_j \ne 0$
leads to acceptance of $H_{0j}$.

\medskip

\noindent The naive $1-\alpha$ confidence interval for $\theta$ is the interval $I(\tilde{K})$.

Suppose that the integer $\ell$ satisfies $q+1 < \ell < p$. Let ${\cal K}^*$ denote the family 
of all subsets of $\{\ell+1, \ldots, p\}$, including the empty set $\emptyset$.
The following theorem paves the way for Theorem 2 which is the main result of this section.

\medskip

\noindent {\bf Theorem 1.} Consider the following 4 cases.

\medskip

\noindent \underbar{Case 1} \ $K^*$ minimizes $AIC(K)$ with respect to $K \in {\cal K}^*$.

\medskip

\noindent \underbar{Case 2} \ $K^*$ minimizes $C_K$ with respect to $K \in {\cal K}^*$.

\medskip

\noindent \underbar{Case 3} \ $K^*$ minimizes $B_K$ with respect to $K \in {\cal K}^*$.

\medskip

\noindent \underbar{Case 4} \ $K^*$ consists of the set of $j \in \{q+1, \ldots, p \}$ for which a $t$-test
of the null hypothesis $H_{0j}: \beta_j=0$ 
against the alternative hypothesis $H_{Aj}: \beta_j \ne 0$
leads to acceptance of $H_{0j}$.

\medskip

\noindent For each of these cases, the coverage probability of the confidence interval
$I(K^*)$ is a function of 
$\frac{1}{\sigma}(\beta_{\ell+1},\ldots, \beta_p)$.

\medskip

\noindent This theorem is proved in Appendix A.

It is intuitively plausible that the wider the class of models that one selects from using
a given model selection procedure, the smaller is the minimum coverage probability of the 
naive $1-\alpha$ confidence interval. The following theorem formalizes this plausible result.
We will use this theorem in Section 3 to derive an easily-computed finite sample upper bound on the minimum coverage
probability of the naive $1-\alpha$ confidence interval.

\medskip

\noindent {\bf Theorem 2.} Consider the following 4 cases.

\medskip

\noindent \underbar{Case 1} \ $\tilde K$ minimizes $AIC(K)$ with respect to $K \in \tilde{\cal K}$. \newline
\phantom{1234567}$K^*$ minimizes $AIC(K)$ with respect to $K \in {\cal K}^*$.

\medskip

\noindent \underbar{Case 2} \ $\tilde K$ minimizes $C_K$ with respect to $K \in \tilde{\cal K}$. \newline
\phantom{1234567}$K^*$ minimizes $C_K$ with respect to $K \in {\cal K}^*$.

\medskip

\noindent \underbar{Case 3} \ $\tilde K$ minimizes $B_K$ with respect to $K \in \tilde{\cal K}$. \newline
\phantom{1234567}$K^*$ minimizes $B_K$ with respect to $K \in {\cal K}^*$.

\medskip

\noindent \underbar{Case 4} \ $\tilde K$ consists of the set of $j \in \{q+1, \ldots, p \}$ for which a $t$-test
of the null hypothesis $H_{0j}: \beta_j=0$ 
against the alternative hypothesis $H_{Aj}: \beta_j \ne 0$
leads to acceptance of $H_{0j}$. 
$K^*$ consists of the set of $j \in \{\ell+1, \ldots, p \}$ for which a $t$-test
of the null hypothesis $H_{0j}: \beta_j=0$ 
against the alternative hypothesis $H_{Aj}: \beta_j \ne 0$
leads to acceptance of $H_{0j}$.

\medskip

\noindent For each of these cases,
the minimum coverage probability of the naive $1-\alpha$ confidence interval
$I(\tilde K)$ is bounded above by the
coverage probability of $I(K^*)$ for each given
$\frac{1}{\sigma}(\beta_{\ell+1},\ldots, \beta_p) \in \mathbb{R}^{p-\ell}$.

\medskip

\noindent This theorem is proved in Appendix B.


\bigskip

\begin{center}
{\large { \bf {3. An easily-computed finite sample upper bound on the 
minimum coverage probability of the naive confidence interval}}}
\end{center}

In this section we present an easily-computed finite sample upper bound on the minimum coverage probability of the 
naive $1-\alpha$ confidence interval. Theorem 2 implies that 
(for each of the methods considered)
this minimum coverage probability is bounded above
by the coverage probability of the 
naive $1-\alpha$ confidence interval for given $\frac{1}{\sigma} \beta_p$ and the model selection procedure applied
only to $\beta_p$. Theorem 3 provides an easily-computed 
expression for the latter coverage probability. This expression is easily minimized numerically with respect to 
$\frac{1}{\sigma} \beta_p$ to obtain the value of the finite sample upper bound on 
the minimum coverage probability of the 
naive $1-\alpha$ confidence interval.

Define the matrix $V$ to be the covariance matrix of $(\hat \Theta, \hat \beta_p)$ divided by $\sigma^2$.
Let $v_{ij}$ denote the $(i,j)$ th element of $V$. Also define the random variable
\begin{equation*}
W = \sqrt{\frac{RSS/(n-p)}{\sigma^2}}
\end{equation*}
and the parameter 
\begin{equation*}
\gamma = \frac{\beta_p}{\sigma\sqrt{v_{22}}}.
\end{equation*}
The random variable $W$ has the same distribution as $\sqrt{Q/(n-p)}$ where $Q \sim \chi^2_{n-p}$.
We have defined $\rho = \text{Corr}(\hat \Theta, \hat \beta_p)$, so that $\rho = v_{12}/\sqrt{v_{11} v_{22}}$.
Define the functions
%
%
%
\begin{align*}
\ell_1(w) &= - t(n-p) w \\
u_1(w) &=  t(n-p) w 
\end{align*}
\begin{align*}
\ell_2(h,w,\rho) &= \rho h - t(n-p+1) \sqrt{\frac{(n-p)w^2+h^2}{n-p+1}} \sqrt{1-\rho^2} \\
u_2(h,w,\rho) &= \rho h + t(n-p+1) \sqrt{\frac{(n-p)w^2+h^2}{n-p+1}} \sqrt{1-\rho^2} 
\end{align*}
Now define the functions
\begin{align*}
k^{\dag}(h,w, \gamma, \rho) &= \Psi \left( \ell_1(w), u_1(w); \rho(h-\gamma), 1-\rho^2 \right )
\\
k(h,w,\gamma, \rho) &= \Psi \left(\ell_2(h,w,\rho), u_2(h,w,\rho); \rho(h-\gamma),1-\rho^2 \right )
\end{align*}
where $\Psi(x, y; \mu, v) = P(x \le Z \le y)$ for $Z \sim N(\mu,v)$. 
Also define \newline
$T = \hat \beta_p/ \big(\sqrt{RSS/(n-p)} \sqrt{v_{22}} \big)$.
We use these definitions in the statement
of the following theorem.

\medskip

\noindent {\bf Theorem 3.} Suppose that ${\cal K}^* = \{ \{p\}, \emptyset\}$.
Consider the following 4 cases.

\medskip

\noindent \underbar{Case 1} \ $K^*$ minimizes $AIC(K)$ with respect to $K \in {\cal K}^*$. Define
\begin{equation*}
d = \sqrt{ \left ( \exp \left ( \frac{2 f(n)}{n} \right ) - 1 \right ) (n-p)}.
\end{equation*}

\medskip

\noindent \underbar{Case 2} \ $K^*$ minimizes $C_K$ with respect to $K \in {\cal K}^*$. Define
$d = \sqrt{2}$.

\medskip

\noindent \underbar{Case 3} \ $K^*$ minimizes $B_K$ with respect to $K \in {\cal K}^*$. Define
$d = 1$.

\medskip

\noindent \underbar{Case 4} \  If $|T| \ge d$ then $K^*=\emptyset$;
otherwise $K^*=\{ p \}$.

\medskip

\noindent In each of these 4 cases, the coverage probability of the confidence interval
$I(K^*)$ is an even function of $\gamma$ and 
is equal to
\begin{equation}
\label{cov_prob_num}
(1-\alpha) 
+ \int_0^{\infty} \int_{-d}^d \big(k(wx,w,\gamma, \rho) - k^{\dag}(wx,w,\gamma, \rho) \big) \, 
\phi(wx - \gamma) \, w \, f_W(w) \, dx \, dw
\end{equation}
where $\phi$ denotes the $N(0,1)$ probability density function and $f_W$ denotes the probability density function 
of $W$. For given $\gamma$, \eqref{cov_prob_num} is an even function of $\rho$.

\medskip

\noindent This theorem is proved in Appendix C. It has the following corollary.

\medskip

\noindent {\bf Corollary 1.} Consider the 4 cases described in Theorem 3. In each of these 4 cases,
the minimum coverage probability of the naive $1-\alpha$ confidence interval is bounded
above by the minimum over $\gamma \ge 0$ of \eqref{cov_prob_num}.

\medskip

That this corollary is a finite sample analogue of Theorem 1 of 
Kabaila \& Leeb (2006) is confirmed as follows. 
The following are conservative model selection procedures: minimizing AIC, minimizing Mallows' $C_P$ and
maximizing adjusted $R^2$. Define $d^{\prime} = \sqrt{2}$ for model selection by minimizing AIC
and by minimizing Mallows' $C_P$. Also define $d^{\prime} = 1$ for model selection by 
maximizing adjusted $R^2$.
Define $z$ by $P(-z \le Z \le z)  = 1 - \alpha$ for $Z \sim N(0,1)$.
Also define $\Delta(a,b) = \Phi(a+b) - \Phi(a-b)$ for all $a, b \in \mathbb{R}$,
where $\Phi$ denotes the $N(0,1)$ distribution function.
Consider $\rho$ and $p$ fixed and $n \rightarrow \infty$. Now $t(n-p) \rightarrow z$ 
as $n \rightarrow \infty$. 
For model selection using AIC, $d \rightarrow d^{\prime}$ as $n \rightarrow \infty$.
It may be shown that, for each of these conservative model selection procedures,
\eqref{cov_prob_num} converges to 
\begin{align}
\label{cov_prob_large_n}
&1-\alpha + \int_{-d^{\prime}}^{d^{\prime}} 
\left ( \Delta \left( \frac{\rho \gamma}{\sqrt{1-\rho^2}},z \right )
- \Delta \left( \frac{\rho (h-\gamma)}{\sqrt{1-\rho^2}}, \frac{z}{\sqrt{1-\rho^2}} \right )
\right ) \phi(h-\gamma) \, dh \notag \\ \notag \\
&= 1-\alpha + 
 \Delta \left( \frac{\rho \gamma}{\sqrt{1-\rho^2}},z \right ) \Delta \left( \gamma, d^{\prime} \right )
- \int_{-d^{\prime}}^{d^{\prime}} \Delta \left( \frac{\rho (h-\gamma)}{\sqrt{1-\rho^2}}, 
\frac{z}{\sqrt{1-\rho^2}} \right ) \phi(h-\gamma) \, dh 
\end{align}
uniformly in $\gamma$ as $n \rightarrow \infty$. 
Now
\begin{equation*}
\int_{-d^{\prime}}^{d^{\prime}} \Delta \left( \frac{\rho (h-\gamma)}{\sqrt{1-\rho^2}}, 
\frac{z}{\sqrt{1-\rho^2}} \right ) \phi(h-\gamma) \, dh 
= P \big(-z \le A \le z, - d^{\prime} \le B \le d^{\prime} \big)
\end{equation*}
where
\begin{equation*}
\left[\begin{matrix} A\\ B \end{matrix}
\right] \sim N \left ( \left[\begin{matrix} 0 \\ \gamma \end{matrix}
\right], \left[\begin{matrix} 1 \quad \rho\\ \rho \quad 1 \end{matrix}
\right] \right ).
\end{equation*}
Define $\tilde A = A + \gamma$ and $\tilde B = B - \gamma$. Thus
\begin{equation*}
\left[\begin{matrix} \tilde B\\ \tilde A \end{matrix}
\right] \sim N \left ( \left[\begin{matrix} 0 \\ \gamma \end{matrix}
\right], \left[\begin{matrix} 1 \quad \rho\\ \rho \quad 1 \end{matrix}
\right] \right )
\end{equation*}
and so 
\begin{align*}
P \big(-z \le A \le z, - d^{\prime} \le B \le d^{\prime} \big)
&= P \big(- d^{\prime} \le \tilde A \le d^{\prime}, -z \le \tilde B \le z,  \big)\\
&= \int_{-z}^{z} \Delta \left( \frac{\gamma + \rho h}{\sqrt{1-\rho^2}}, 
\frac{d^{\prime}}{\sqrt{1-\rho^2}} \right ) \phi(h) \, dh 
\end{align*}
Thus \eqref{cov_prob_large_n} is equal to (4) of Kabaila \& Leeb (2006).
This shows that the finite sample upper bound stated in Corollary 1 converges to the 
large sample upper bound
(3) of Kabaila \& Leeb (2006) as $n \rightarrow \infty$.

The following result provides an explicit formula for the upper bound described in Corollary 1 for the
particular case that $\rho=1$. The proof of this result is omitted for the sake of brevity.

\medskip

\noindent {\bf Theorem 4.} Suppose that $\rho=1$. Let $d$ be as defined in the statement of Theorem 3.
The upper bound, described in Corollary 1, to the 
minimum coverage probability of the naive $1-\alpha$ confidence interval is 
\begin{equation*}
2 \int_0^{\infty} \big( \Phi(t(n-p) w) - \Phi(d w) \big) f_W(w) dw
\end{equation*}
when $d < t(n-p)$ and is 0 when $d \ge t(n-p)$.

\bigskip

\begin{center}
{\large { \bf {4. Numerical illustrations}}}
\end{center}

The integrand of the double integral in \eqref{cov_prob_num} is a smooth function of $(x,w)$ and
so it is easily computed numerically. Let $m=n-p$
and remember that $\rho=\text{Corr}(\hat \Theta, \hat \beta_p)$,
where $p$ maximizes $|\text{Corr}(\hat \Theta, \hat \beta_j)|$
with respect to $j \in \{ q+1, \ldots, p\}$. 
For given $p$, $m$, $\alpha$ and $\rho$, we minimize 
\eqref{cov_prob_num} numerically with respect to $\gamma \ge 0$ to obtain the upper bound (described in 
Corollary 1) to the 
minimum coverage probability of the naive $1-\alpha$ confidence interval $I(\tilde K)$. 
The following are conservative model selection procedures: minimizing Mallows' $C_P$, maximizing adjusted $R^2$ 
and minimizing AIC. For the numerical illustrations for these procedures 
described in this section we include the case $m=\infty$. For this case, we use the 
large sample upper bound to the minimum coverage probability of the naive
$1-\alpha$ confidence interval derived by Kabaila \& Leeb (2006).
Programs for computing these upper bounds have been written in MATLAB (including the use of
the Optimization and Statistics toolboxes).

For model selection by minimizing Mallows' $C_P$ or maximizing adjusted $R^2$, $d$ is a fixed number
that does not depend on either $p$ or $m$. In this case, the upper bound (described in Corollary 1)
to the minimum coverage
probability of the naive $1-\alpha$ confidence interval is, for given $|\rho|$,
a function of $m$. 
Plots of this upper bound as a function of $|\rho|$, 
for model selection by minimizing Mallows' $C_P$ and by maximizing adjusted $R^2$,
were prepared for $\alpha \in \{0.1, 0.05, 0.02\}$
and $m=1, 2, 3, 4, 5, 10, 20, 50, 1000$ and $\infty$.
For each value of $\alpha$ and $m$ considered,
this upper bound was found to be a continuous decreasing function of $|\rho|$
that is far below $1-\alpha$ when $|\rho|$ is close to 1.
This finding is illustrated by Figures 1 and 2.
Figure 1 is a plot of this upper bound as a function of $|\rho|$ 
for model selection by minimizing Mallows' $C_P$ and for $m = 5, 20, 50, 1000$ and $\infty$. 
Figure 2 is a plot of this upper bound as a function of $|\rho|$ 
for model selection by maximizing adjusted $R^2$ and for $m = 5, 20, 50, 1000$ and $\infty$.

Now consider model selection using AIC. When $n$ is large and $p$ is small compared to $n$, 
$d$ is approximately equal to $\sqrt{2}$ and the upper bound described by Corollary 1 
is approximately equal to this upper bound for model selection by minimizing Mallows' $C_P$.
Plots of this upper bound as a function of $|\rho|$, for model selection by minimizing AIC,
were prepared for $\alpha=0.05$, $p \in \{2,3,4,7,10\}$ and $m=1, 2, 3, 4, 5, 10, 20, 50, 
1000$ and $\infty$. 
For each value of $p$ and $m$ considered,
this upper bound was found to be a continuous decreasing function of $|\rho|$
that is far below $1-\alpha$ when $|\rho|$ is close to 1.
This finding is illustrated by Figure 3.
This figure is a plot of the upper bound described by Corollary 1 as a function of $|\rho|$
for model selection by minimizing AIC, for $\alpha=0.05$, $p=10$ and $m = 5, 20, 50, 1000$ and $\infty$. 
For the real life data example considered by Kabaila \& Leeb (2006, section 3), $p=10$ and
$m=20$.

Finally, consider model selection using BIC. Since this model selection procedure is consistent,
the large sample upper bound to the minimum coverage probability of the naive
$1-\alpha$ confidence interval, derived by Kabaila \& Leeb (2006), does not apply.
Plots of this upper bound as a function of $|\rho|$, for model selection by minimizing BIC,
were prepared for $\alpha=0.05$, $p \in \{2,3,4,7,10\}$ and $m=1, 2, 3, 4, 5, 10, 20, 50, 
1000$ and $10,000$. 
For each value of $p$ and $m$ considered,
this upper bound was found to be a continuous decreasing function of $|\rho|$
that is far below $1-\alpha$ when $|\rho|$ is close to 1.
This finding is illustrated by Figure 4.
This figure is a plot of the upper bound described by Corollary 1 as a function of $|\rho|$
for model selection by minimizing BIC, for $\alpha=0.05$, $p=10$ and $m = 5, 20, 50, 1000$
and $10,000$.

\medskip


\FloatBarrier

\begin{figure}[t]
\label{Figure1}
    \centering
    \includegraphics[scale=0.65]{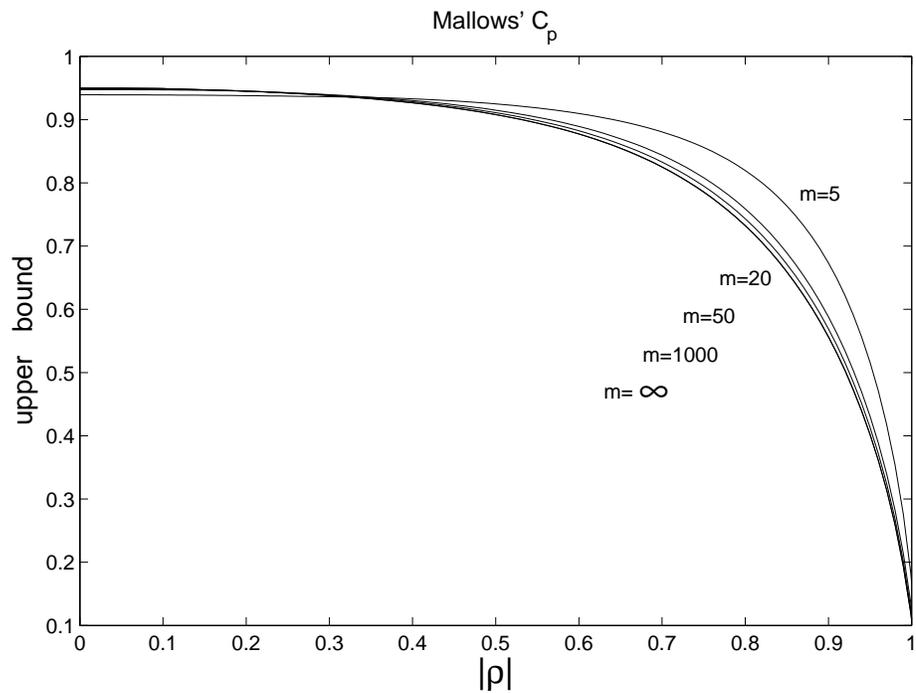}
    \caption{Plot of the upper bound, stated in Corollary 1, 
    on the coverage probability of the naive 95\% confidence 
    interval against $|\rho|$ when model selection is by minimization of
    Mallows' $C_P$. 
    Here  $m = n-p = 5, 20, 50, 1000$ and
    $\infty$.}
\end{figure}
%
%
\begin{figure}[h]
\label{Figure2}
    \centering
    \includegraphics[scale=0.65]{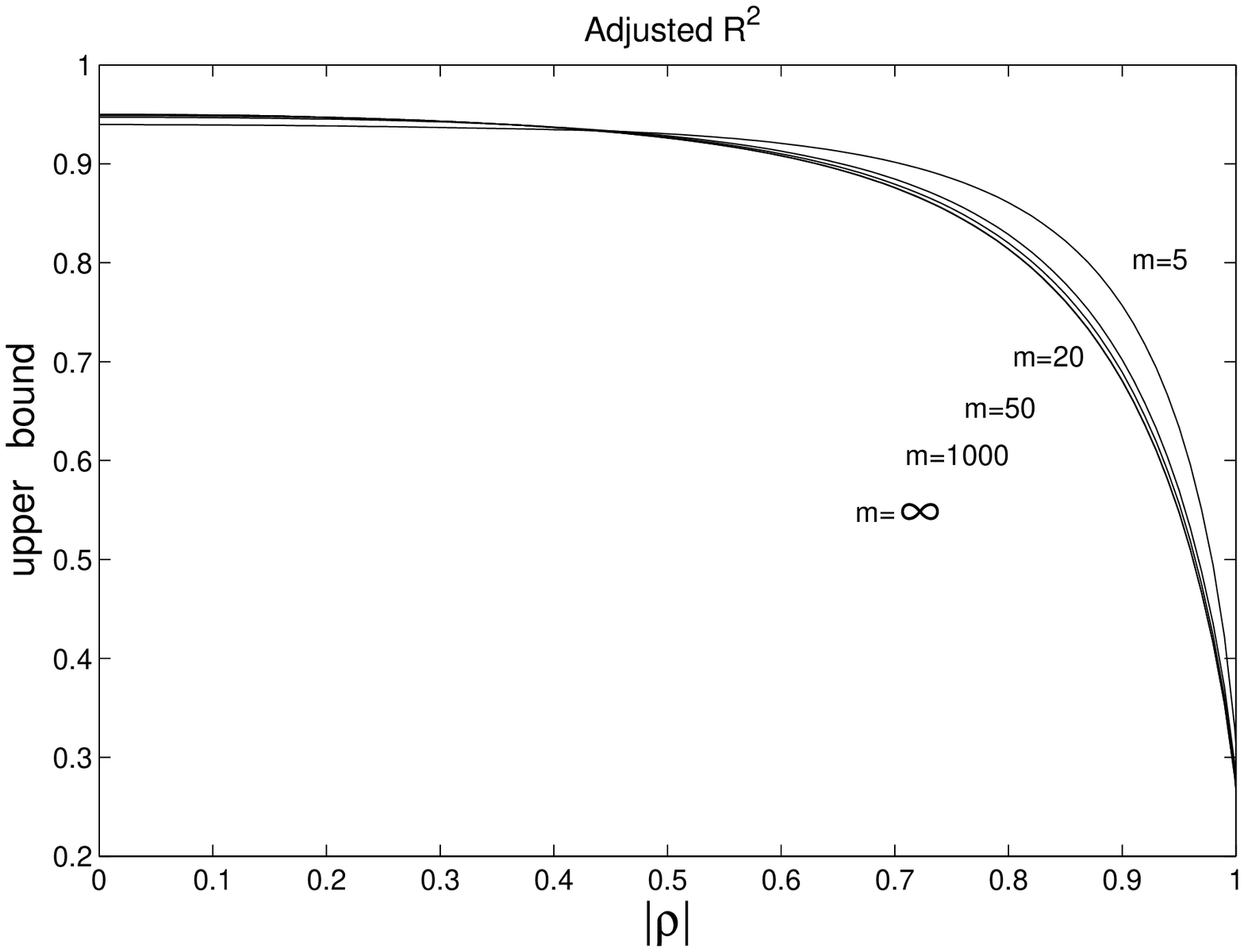}
    \caption{Plot of the upper bound, stated in Corollary 1, 
    on the coverage probability of the naive 95\% confidence 
    interval against $|\rho|$ when model selection is by maximization of
    adjusted $R^2$. 
    Here  $m = n-p = 5, 20, 50, 1000$ and
    $\infty$.}
\end{figure}
\newpage
\FloatBarrier
\begin{figure}[t]
\label{Figure3}
    \centering
    \includegraphics[scale=0.64]{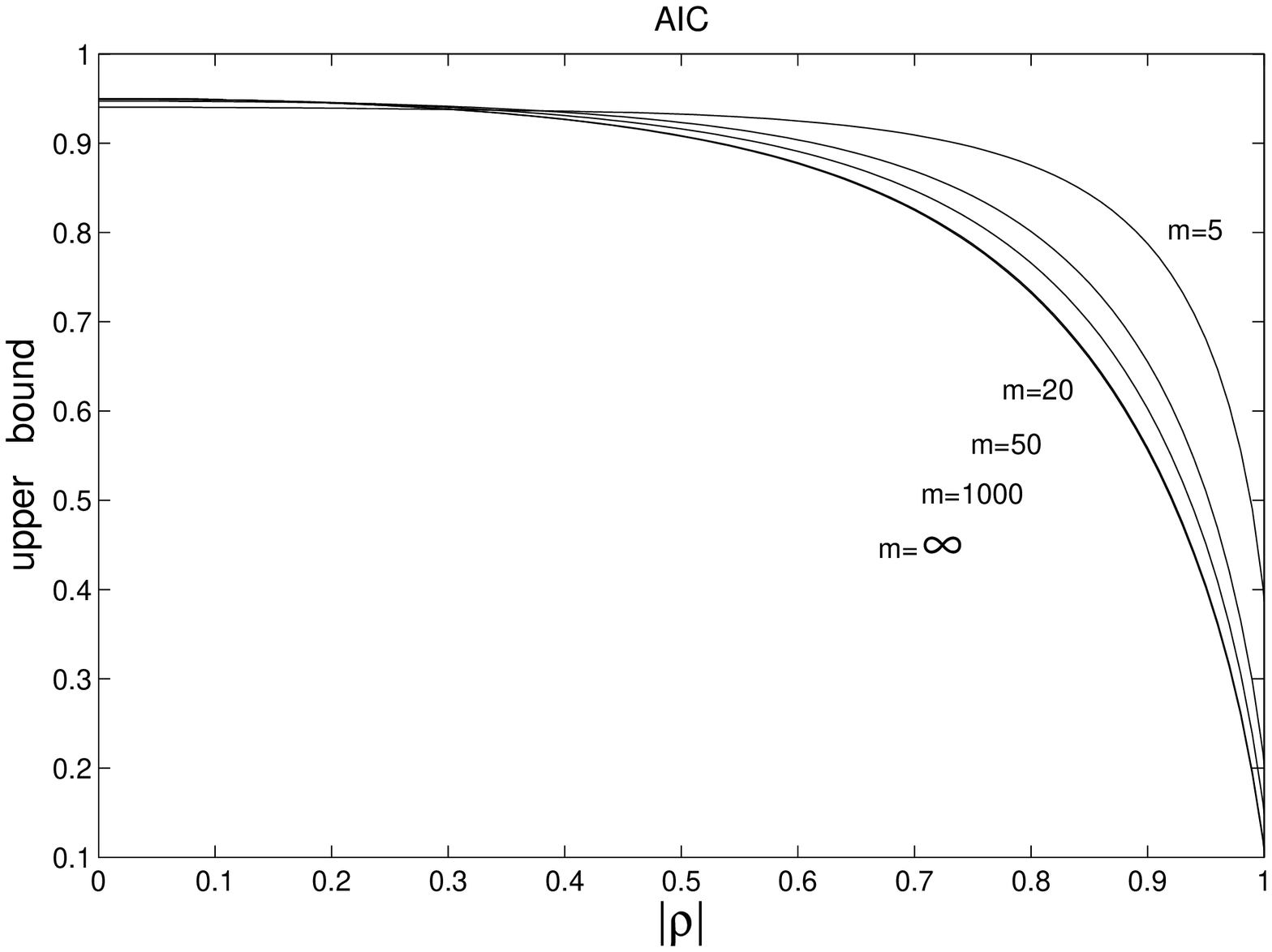}
   \caption{Plot of the upper bound, stated in Corollary 1, 
    on the coverage probability of the naive 95\% confidence 
    interval against $|\rho|$ when model selection is by minimization of
    AIC. 
    Here  $p=10$ and $m = n-p = 5, 20, 50, 1000$ and
    $\infty$.}
\end{figure}
%
\begin{figure}[h]
\label{Figure4}
    \centering
    \includegraphics[scale=0.64]{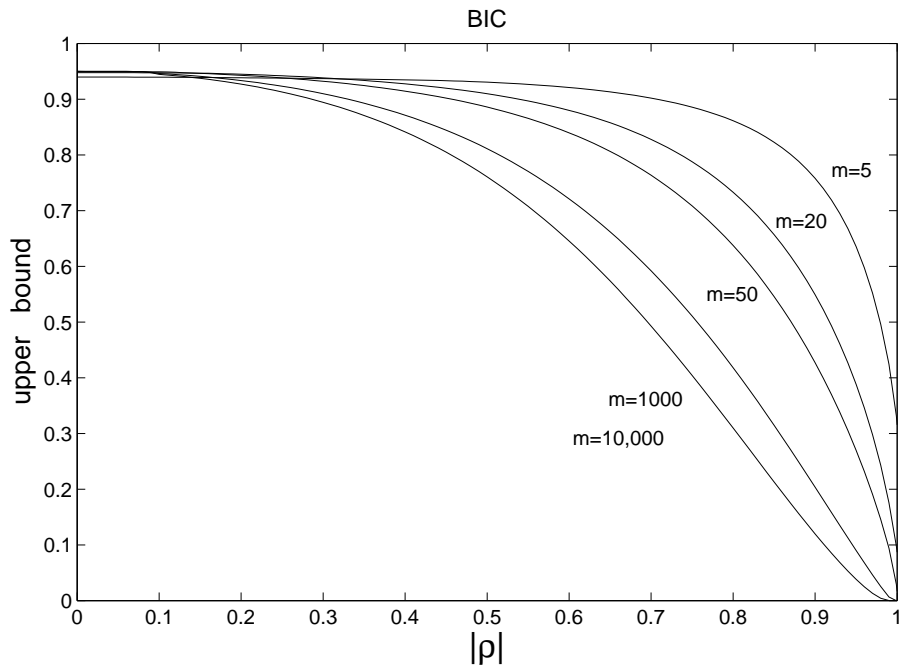}
   \caption{Plot of the upper bound, stated in Corollary 1, 
    on the coverage probability of the naive 95\% confidence 
    interval against $|\rho|$ when model selection is by minimization of
    BIC. 
    Here  $p=10$ and $m = n-p = 5, 20, 50, 1000$ and $10,000$.}
\end{figure}

\FloatBarrier


\begin{center}
{\large { \bf {5. Conclusion}}}
\end{center}

For a given design matrix $X$ and a wide variety of model selection procedures, 
the efficient Monte Carlo simulation methods of 
Kabaila (2005) 
and Giri \& Kabaila (2007) 
provide valuable information about the minimum coverage probability 
of the naive $1-\alpha$ confidence interval. What is also of interest, however, is to delineate general categories
of design matrices $X$ and model selection procedures for which this confidence interval has poor coverage properties.
The first such delineation, for the complicated kinds of model selection procedures used in practice, 
results from the upper bound on the minimum coverage probability of this confidence interval due to  
Kabaila \& Leeb (2006). This upper bound, however, is valid only in large samples and 
applies only to conservative model selection
procedures. The present paper presents a finite sample analogue of this upper bound that is applicable to 
a wide variety of model selection procedures and provides a delineation that is valid for finite samples.

\bigskip

\begin{center}
{\large { \bf {Appendix A: Proof of Theorem 1}}}
\end{center}

In this appendix we prove Theorem 1. The proof is in 2 parts.

\smallskip

\noindent \underbar{Part 1} \ For each of the Cases 1--4, $K^*$ is determined by the following set of
random variables
\begin{equation*}
\left\{ \frac{RSS}{\sigma^2} \right \} \cup
\left\{ \frac{RSS_K}{\sigma^2}: K \in {\cal K}^* \right \}.
\end{equation*}
By Theorem 1(c) and the proof of Theorem
1(e) of Kabaila (2005), in each of these cases, $K^*$ is determined by
\begin{equation*}
\left ( \frac{RSS}{\sigma^2}, \eta^{''}, \frac{1}{\sigma} (\beta_{\ell+1},\ldots, \beta_p) \right )
\end{equation*}
where the random vector $\eta^{''}$ is defined by Kabaila (2005, p. 552).

\smallskip

\noindent \underbar{Part 2} \ It follows from Part 1 and the proof of Theorem 1(f) of Kabaila (2005) that, 
in each of the 4 cases, $P \big( \theta \in I(K^*) \big)$ is a function of 
$\frac{1}{\sigma} (\beta_{\ell+1},\ldots, \beta_p)$.

\newpage

\begin{center}
{\large { \bf {Appendix B: Proof of Theorem 2}}}
\end{center}

In this appendix we prove Theorem 2. The proof is in 2 parts.

\smallskip

\noindent \underbar{Part 1}\ \ Suppose that $K \ne \emptyset$.
It is well-known (see e.g. Graybill (1976, p.222)) that
\begin{equation*}
RSS_K=RSS+(H_K\hat{\beta})^T \big(H_K(X^TX)^{-1}H^T_K \big)^{-1}
H_K\hat{\beta}.
\end{equation*}
Thus
\begin{equation*}
\frac{RSS_K}{\sigma^2}=\frac{RSS}{\sigma^2}+V_K.
\end{equation*}
where
\begin{equation*}
V_K=(H_K \textstyle{\frac{1}{\sigma}} \hat{\beta})^T \big(H_K(X^TX)^{-1}H^T_K \big)^{-1}
H_K \textstyle{\frac{1}{\sigma}} \hat{\beta}.
\end{equation*}
By a well-known result (see e.g. Graybill (1976, p.127)), $V_K$ has a noncentral
chi-squared distribution with degrees of freedom $|K|$ and noncentrality
parameter $\lambda = \frac{1}{2}(H_K \textstyle{\frac{1}{\sigma}} {\beta})^T \big(H_K(X^TX)^{-1}H^T_K \big)^{-1}
H_K \textstyle{\frac{1}{\sigma}} {\beta}$ (in the notation for noncentral 
chi-squared distributions used by Graybill (1976)). 

In Cases 1--3, express $\tilde K$ in terms of the following set of random variables
\begin{equation*}
\left\{ \frac{RSS}{\sigma^2} \right \} \cup
\left\{ V_K: K \in \tilde{\cal K} \right \}.
\end{equation*}
In Case 4, 
express $\tilde K$ in terms of the following set of random variables
\begin{equation*}
\left\{ \frac{RSS}{\sigma^2} \right \} \cup
\left\{ V_{\{ j \}}: j \in \{q+1, \ldots, p \} \right \}.
\end{equation*}
Note that $RSS/\sigma^2$ and $V_K$ are independent random variables and $RSS/\sigma^2 \sim \chi^2_{n-p}$.

\smallskip

\noindent \underbar{Part 2}\ \ Fix
$\frac{1}{\sigma} (\beta_{\ell+1},\ldots, \frac{1}{\sigma} \beta_p)$.
Choose
$|\frac{1}{\sigma} \beta_{q+1}|= \cdots = |\frac{1}{\sigma} \beta_{\ell}|$
and consider 
$|\frac{1}{\sigma} \beta_{q+1}|= \cdots = |\frac{1}{\sigma} \beta_{\ell}| \rightarrow \infty$.
Define ${\cal J}$ to be the family of sets that belong to $\tilde{\cal K}$
and include at least one element of $\{q+1, \ldots, \ell\}$.

\noindent (a) Using the expression for $\tilde K$ found in Part 1,
it may be shown that for each of the 4 cases and for each $K \in {\cal J}$,
\begin{equation*}
P(\tilde K =K) \rightarrow 0
\end{equation*}
as $|\frac{1}{\sigma} \beta_{q+1}|= \cdots = |\frac{1}{\sigma} \beta_{\ell}| \rightarrow \infty$.
For example, for Case 1 
minimizing $AIC(K)$ with respect to $K \in \tilde{\cal K}$ is equivalent to minimizing
\begin{equation*}
IC(K) = n \ln \left(\frac{RSS}{\sigma^2} + V_K\right)
+2(p-\vert K\vert )f(n)
\end{equation*}
with respect to $K \in \tilde{\cal K}$.
Thus, for each $K \in {\cal J}$,
\begin{equation*}
P(\tilde K =K) \le P\big(IC(K) \le IC(\emptyset)\big) \rightarrow 0
\end{equation*}
as $|\frac{1}{\sigma} \beta_{q+1}|= \cdots = |\frac{1}{\sigma} \beta_{\ell}| \rightarrow \infty$.
Hence, in each of the 4 cases, $P(\tilde K \in {\cal J}) \rightarrow 0$ as 
$|\frac{1}{\sigma} \beta_{q+1}|= \cdots = |\frac{1}{\sigma} \beta_{\ell}| \rightarrow \infty$.

\medskip

\noindent (b) Observe that the minimum value of $P(\theta \in I(\tilde K))$ is bounded above by
\begin{align*}
P(\theta \in I(\tilde K)) &= P \left ( \cup_{K \in {\cal J}^c} \left ( \{ \theta \in I(K) \} \cap \{\tilde K = K \} \right) \right ) 
+ P \left ( \cup_{K \in {\cal J}} \left ( \{ \theta \in I(K) \} \cap \{\tilde K = K \} \right) \right )\\
&\le P \left ( \cup_{K \in {\cal J}^c} \left ( \{ \theta \in I(K) \} \cap \{\tilde K = K \} \right) \right )
+ P(\tilde K \in {\cal J})\\
&\le P \left ( \cup_{K \in {\cal J}^c} \left ( \{ \theta \in I(K) \} \cap \{K^* = K \} \right) \right )
+ P(\tilde K \in {\cal J})\\
&= P (\theta \in I(K^*) ) + P(\tilde K \in {\cal J})
\end{align*}
since $\{ \tilde K = K \} \subset \{ K^* = K \}$ for each $K \in {\cal J}^c$.
By choosing $|\frac{1}{\sigma} \beta_{q+1}|= \cdots = |\frac{1}{\sigma} \beta_{\ell}| \rightarrow \infty$,
we see that the minimum value of $P(\theta \in I(\tilde K))$ is bounded above by $P (\theta \in I(K^*) )$.
 
\bigskip

\begin{center}
{\large { \bf {Appendix C: Proof of Theorem 3}}}
\end{center}

In this appendix we prove Theorem 3. Define the random variables 
\begin{equation*}
G = \frac{\hat \Theta - \theta}{\sigma \sqrt{v_{11}}} \qquad \text{and} \qquad
H = \frac{\hat \beta_p}{\sigma \sqrt{v_{22}}}.
\end{equation*}
Note that $T=H/W$.
%
%
In each of the 4 cases, 
$K^* = \emptyset$ if $|T| \ge d$ and $K^* = \{ p \}$ otherwise.
It is straightforward to show that the confidence interval $I(K^*)$ for $\theta$ is
\begin{equation*}
\left [ \hat \Theta - t(n-p) \sqrt{v_{11}} \sqrt{\frac{RSS}{n-p}}, \,
\hat \Theta + t(n-p) \sqrt{v_{11}} \sqrt{\frac{RSS}{n-p}} \right ]
\end{equation*}
if $|T| \ge d$ and 
\begin{align*}
\Bigg [ \hat \Theta -  \frac{v_{12}}{v_{22}} \hat \beta_p 
- &t(n-p+1) \sqrt{\frac{RSS + (\hat \beta_p^2 / v_{22})}{n-p+1}}\sqrt{v_{11} - \frac{v_{12}^2}{v_{22}}}, \\
&\hat \Theta -  \frac{v_{12}}{v_{22}} \hat \beta_p 
+ t(n-p+1) \sqrt{\frac{RSS + (\hat \beta_p^2 / v_{22})}{n-p+1}}\sqrt{v_{11} - \frac{v_{12}^2}{v_{22}}} \, \Bigg ]
\end{align*}
otherwise. Note that 
\begin{equation}
\label{model_G_H}
\tag{C.1}
\left[\begin{matrix} G\\ H \end{matrix}
\right] \sim N \left ( \left[\begin{matrix} 0 \\ \gamma \end{matrix}
\right], \left[\begin{matrix} 1 \quad \rho\\ \rho \quad 1 \end{matrix}
\right] \right ).
\end{equation}
It may be shown that the coverage probability of $I(K^*)$ is equal to 
\begin{align}
\label{cov_first}
P &\left ( \{ \ell_1(W) \le G \le u_1(W) \} \cap \left \{ \frac{|H|}{W} \ge d \right \} \right ) \notag \\
&+ P \left ( \{ \ell_2(H,W,\rho) \le G \le u_2(H,W,\rho) \} \cap \left \{ \frac{|H|}{W} < d \right \} \right ).
\tag{C.2}
\end{align}
Remember that $\ell_1$, $u_1$, $\ell_2$ and $u_2$ are defined at the start of Section 3.
Using the fact that 
\begin{equation*}
\left[\begin{matrix} -G\\ -H \end{matrix}
\right] \sim N \left ( \left[\begin{matrix} 0 \\ -\gamma \end{matrix}
\right], \left[\begin{matrix} 1 \quad \rho\\ \rho \quad 1 \end{matrix}
\right] \right ).
\end{equation*}
it may be shown that \eqref{cov_first} is an even function of $\gamma$.

The random vectors $(G,H)$ and $W$ are independent. It follows from \eqref{model_G_H}
that the probability density function of $H$, evaluated at $h$, is $\phi(h-\gamma)$. Thus
\begin{align}
\label{first_term_cov}
 P &\left ( \{ \ell_1(W) \le G \le u_1(W) \} \cap \left \{ \frac{|H|}{W} \ge d \right \} \right ) \notag \\
 &= \int_0^{\infty} \int_{\{|h| > dw \}} \int_{\ell_1(w)}^{u_1(w)}
f_{G|H}(g|h) \, dg \, \phi(h-\gamma)\,  dh \, f_W(w) \, dw
\tag{C.3}
\end{align}
where $f_{G|H}(g|h)$ denotes the probability density function of $G$ conditional on $H=h$, evaluated at $g$.
By \eqref{model_G_H}, the probability distribution of $G$ conditional on $H=h$ is
$N \big( \rho (h - \gamma), 1 - \rho^2 \big )$. It follows that \eqref{first_term_cov} is equal
\begin{equation}
\label{first_term_cov_simpler}
  \int_0^{\infty} \int_{\{|h| > dw \}} k^{\dag}(h,w, \gamma, \rho)
\, \phi(h-\gamma)\, f_W(w) \,  dh  \, dw.
\tag{C.4}
\end{equation}
The standard $1-\alpha$ confidence interval $I(\emptyset)$ for $\theta$ has coverage probability $1-\alpha$,
so that $1 - \alpha = P \big(\ell_1(W) \le G \le u_1(W) \big)$. Thus
\eqref{first_term_cov_simpler} is equal to
\begin{equation*}
  (1-\alpha) -\int_0^{\infty} \int_{-dw}^{dw} k^{\dag}(h,w,\gamma, \rho)
\, \phi(h-\gamma)\, f_W(w) \,  dh  \, dw.
\end{equation*}
Similarly,
\begin{align*}
P &\left ( \{ \ell_2(H,W,\rho) \le G \le u_2(H,W,\rho) \} \cap \left \{ \frac{|H|}{W} < d \right \} \right )\\
&= \int_0^{\infty} \int_{-dw}^{dw} k(h,w,\gamma, \rho)
\, \phi(h-\gamma)\, f_W(w) \,  dh  \, dw.
\end{align*}
Hence, $P(\theta \in I(\bar K))$ is equal to 
\begin{equation*}
  (1-\alpha) + \int_0^{\infty} \int_{-dw}^{dw} \big (k(h,w,\gamma, \rho) - k^{\dag}(h,w,\gamma, \rho) \big )
\, \phi(h-\gamma)\, f_W(w) \,  dh  \, dw.
\end{equation*}
The result follows by changing the variable of integration in the inner integral from $h$ to
$x=h/w$.

\noindent That, for given $\gamma$, \eqref{cov_prob_num} is an even function of $\rho$ follows from the
fact that $\Phi(b) - \Phi(a) = \Phi(-a) - \Phi(-b)$ for all $a, b \in \mathbb{R}$.

\bigskip

\begin{center} {\large {\sl References}} \end{center}

\rf ARABATZIS, A.A., GREGOIRE, T.G., \& REYNOLDS, M.R. (1989). 
Conditional estimation of the mean following rejection of a two sided test.
{\it Communications in Statistics - Theory and Methods} {\bf 18}, 4359--4373.

\rf BREIMAN, L. (1992). The little bootstrap and other methods for dimensionality 
selection in regression: X-fixed prediction error. {\it Journal of the American Statistical Association}
{\bf 87}, 738--754.

\rf CHIOU, P. (1997). Interval estimation of scale parameters following a pre-test for two
exponential distributions. {\it Computational Statistics \& Data Analysis} {\bf 23}, 477--489.

\rf CHIOU, P., \& HAN, C-P. (1995a). Conditional interval estimation of the exponential
location parameter following rejection of a pre-test. {\it Communications in Statistics - Theory and Methods} 
{\bf 24}, 1481--1492.

\rf CHIOU, P., \& HAN, C-P. (1995b). Interval estimation of error variance
following a preliminary test in one-way random model. 
{\it Communications in Statistics - Simulation and Computation} 
{\bf 24}, 817--824.

\rf GIRI, K., \& KABAILA, P. (2007). The Coverage Probability of Confidence Intervals in $2^r$ Factorial
Experiments After Preliminary Hypothesis Testing. 
To appear in {\it Australian \& New 
Zealand Journal of Statistics}.

\rf GRAYBILL, F. A. (1976). {\it Theory and Application of the Linear Model}.
Pacific Grove CA: Duxbury.

\newpage

\rf HAN, C-P. (1998). Conditional confidence intervals of regression coefficients
following rejection of a preliminary test. In {\it Applied Statistical Science}, Vol III, eds. S.E. Ahmed,
M. Ahsanullah, and B.K. Sinha, papers in Honours of A.K.Md.E. Saleh: Nova Science, pp. 193--202.

\rf HURVICH, C.M., \& TSAI, C-L. (1990). The impact of model selection on inference in linear regression,''
{\it The American Statistician} {\bf 44}, 214--217.


\rf KABAILA, P. (1995). The effect of model selection on confidence regions and prediction regions. 
{\it Econometric Theory} {\bf 11}, 537--549.

\rf KABAILA, P. (1998). Valid confidence intervals in
regression after variable selection. {\it Econometric Theory} {\bf 14},
463--482.

\rf KABAILA, P. (2005). On the coverage probability of confidence intervals
in regression after variable selection. {\it Australian \& New Zealand Journal of Statistics}
{\bf 47}, 549--562.

\rf KABAILA, P., \& LEEB, H. (2006). On the large-sample minimum coverage 
probability of confidence intervals after
model selection. {\it Journal of the American Statistical Association} {\bf 101}, 619--629.

\rf LEEB, H., \& P\"OTSCHER, B. M. (2005). Model selection and inference: facts and fiction. 
{\it Econometric Theory} {\bf 21}, 21--59.

\rf REGAL, R.R., \& HOOK, E.B. (1991). The effects of model selection on confidence intervals
for the size of a closed population. {\it Statistics in Medicine} {\bf 10}, 717--721.

\end{document}